\documentclass[10pt,twoside]{amsart}
\usepackage{amssymb,graphicx,verbatim}
\newtheorem{thm}{Theorem}
\newtheorem{lem}[thm]{Lemma}
\newtheorem{prop}[thm]{Proposition}

\newtheorem{defi}[thm]{Definition}
\newcommand{\R}{\mathbb{R}}

\newcommand{\C}{\mathbb{C}}

\newcommand{\finprf}{\unskip\null\hfill$\square$\vskip 0.3cm}

\newcommand{\be}{\begin{equation}}
\newcommand{\ee}{\end{equation}}

\begin{document}


\title[Self-adjointness for Dirac operators via Hardy-Dirac inequalities]{Self-adjointness of Dirac operators via Hardy-Dirac  inequalities}

\author[M.J. Esteban]{Maria J. Esteban$^1$}
\address{$^1$Ceremade, Universit\'e Paris Dauphine, Place de Lattre de Tassigny, F-75775 Paris C\'edex 16, France}
\email{esteban@ceremade.dauphine.fr}

\author[M. Loss]{Michael Loss$^2$}
\address{$^2$School of Mathematics, Georgia Institute of Technology, Atlanta, GA 30332, USA}
\email{loss@math.gatech.edu}

\begin{abstract}
Distinguished selfadjoint extension of Dirac operators are constructed
for a class of potentials including Coulombic ones up to the critical case, $-|x|^{-1}$. The method uses Hardy-Dirac inequalities and quadratic form techniques.
\end{abstract}
\subjclass[2000]{}
\date{\today}
\maketitle
\thispagestyle{empty}


\section{Introduction.}

In  the units in which both the speed of light $c$ and  Planck's constant  $\hbar$ are equal to $1$,
the Dirac operator in the presence of an external electrostatic potential $V$   is given by
\be\label{one}
H_0+V \quad\mbox{with}\quad H_0:=-i\, \alpha \cdot\nabla+\beta \; .
\ee

In (\ref{one}), $\alpha_1$, $\alpha_2$, $\alpha_3$ and $\beta$ are $4 \times 4$ complex matrices, whose standard form (in $2\times 2$ blocks) is
$$
\beta=\left( \begin{matrix} {
\mathbb I} & 0 \\ 0 & -{\mathbb I} \\ \end{matrix} \right),\quad\alpha_k=\left( \begin{matrix} 0 &\sigma_k \\ \sigma_k &0 \\ \end{matrix}\right) \qquad (k=1, 2, 3)\;,
$$
where ${
\mathbb I}=\left( \begin{matrix} 1 & 0 \\ 0 & 1 \end{matrix} \right)$ and $\sigma_k$ are the Pauli matrices:
$$\sigma _1=\left( \begin{matrix} 0 & 1 \\ 1 & 0 \\ \end{matrix} \right),\quad \sigma_2=\left( \begin{matrix} 0 & -i \\ i & 0 \\ \end{matrix}\right),\quad \sigma_3=\left( \begin{matrix} 1 & 0\\ 0 &-1\\ \end{matrix}\right) .$$

If $V$ is a bounded function which tends to $0$ at infinity, one can easily prove that the operator $H_0+V$ with domain $H^1(\R^3, \C^4)$ is self-adjoint. If $V$ has singularities, as it is the case for instance  in  many models for atoms, one is interested in defining self-adjoint extensions of $T:= H_0+V_|{_{C^\infty_0(\R^3,\C^4)}}$. The method used to do this depends on the singularity. Let us for instance consider Coulomb potentials $-\nu/|x|$, $\nu>0$. Then for $\nu\in (0, \pi/2]$  one can use the pseudo-Friedrich extension method to define an extension which satisfies 
\be\label{domain}{\mathcal D}(H_0+V)\subset {\mathcal D}(|H_0|^{1/2})=H^{1/2}(\R^3, \C^4)\,.
\ee
 This result is obtained by using Kato's inequality :
$$|H_0|\geq \frac2{\pi |x|}\,.$$
Actually one can prove that $H_0-\frac\nu{|x|}$ defined on $C^\infty_0(\R^3,\C^4)$ is essentially self-adjoint if $\nu<\sqrt{3}/2$ (\cite{Schmincke-72B}). 

When the singularities are stronger, that is, if $\nu\geq \sqrt{3}/2$,  other methods need to be used. Various works have dealt with this issue, and it appears that for potentials $V$ which have only a singularity at the origin, the condition 
\begin{equation}\label{sotto1}
\,\sup_{x\ne 0}\, |x|\,|V(x)|<1\,,
\end{equation}
 is sufficient to define a distinguised selfadjoint extension of $H_0+V_|{_{C^\infty_0(\R^3\setminus\{0\}, \C^4)}}$. This has been done by means of cut-off methods by W\"ust \cite{Wust-73, Wust-75, Wust-77} in the case of semibounded potentials $V$ and by Schmincke \cite{Schmincke-72} without the assumption of semiboundedness. These extensions are characterized by the fact that the domain is contained in ${D(T^*)\cap D(r^{-1/2})}$. On the other hand, under the same assumption \eqref{sotto1} Nenciu proved in \cite{Nenciu-76} the existence of a unique selfadjoint extension $\tilde T$ with domain contained in $H^{1/2}(\R^3, \C^4)$. Finally, in the case of semibounded potentials satisfying \eqref{sotto1}, Klaus and W\"ust proved in  \cite{Klaus-Wust-78} that all the aforementioned self-adjoint extensions coincide. 
 
 In this paper we address the same question from a different perspective, and we prove that the existence of a distinguished self-adjoint extension of Dirac operators 
$H_0+V_|{_{C^\infty_0(\R^3, \C^4)}}$ can be achieved whenever there exists a Hardy-like inequality involving the operators $H_0$ and $V$. This is somewhat surprising since the Dirac operator is not semibounded. This should be contrasted with the case $-\Delta -\frac\mu{|x|^2}$ where the classical Hardy inequality defines the threshold of self-adjointness. 

In \cite{Dolbeault-Esteban-Sere-00B}  Hardy-like inequalities for Dirac operators were proved. They allow potentials which cover all admissible Coulomb like singularities, that is, all the potentials $V=\frac\nu{|x|}\,, \nu\in (0, 1]$.
As a result, our method allows to overcome the limitation contained in assumption
\eqref{sotto1} in a natural way by constructing distinguished
selfadjoint extensions for potentials satisfying
\begin{equation}\label{sotto2}
\,\sup_{x\ne 0}\, |x|\,V(x) \le 1\,,
\end{equation}
which is the maximal possible range.

One example of a Hardy-Dirac inequality is given by the following theorem.
\begin{thm}[\cite{Dolbeault-Esteban-Sere-00B}]\label{TR5} {\rm} Let $V$ be a  function satisfying 
\be\label{V1-V2} 
\lim_{|x|\to +\infty}V(x)=0\quad\mbox{and}\quad -\frac{\nu}{|x|}-c_1\leq V\leq \Gamma=\sup(V)\;,
\ee
with $\nu\in (0,1)$, $c_1$, $\Gamma \in\R$, $c_1$, $\Gamma \geq 0$, $c_1+\Gamma-1< \sqrt{1-\nu^2}$.
There exists a constant $c(V) \in (-1,1)$ such that
for all $\varphi \in C^\infty_c (\R^3, \C^2)$,
\be\label{R5} 
\int_{\R^3}\left(\frac{|{\boldsymbol{\sigma}}\cdot{\boldsymbol{\nabla}}\varphi|^2}{1+c(V)-V}\ +\ \big(1-c(V)+V\big)\,|\varphi|^2\right)\,dx\ \geq\ 0 \;.
\ee
\end{thm}

Note that the argument in \cite{Dolbeault-Esteban-Sere-00B} proceeded in an indirect way. Using spectral analysis of
the Dirac operator, i.e., its selfadjointness, Theorem \ref{TR5} was established for the case $\nu \in (0,1)$. From this, a simple limiting argument yields the Hardy-Dirac inequality for the case $\nu=1$,
\be
\label{HD} 
\int_{\R^3}\left(\frac{|{\boldsymbol{\sigma}}\cdot{\boldsymbol{\nabla}}\varphi|^2}{1+\frac{1}{|x|}}\ +\ \,|\varphi|^2\right)\,dx\ \geq\ \int_{\R^3} \frac{|\varphi|^2}{|x|} \,dx \;.
\ee
This functional inequality was later derived in \cite{Dolbeault-Esteban-Loss-Vega-04} directly, i.e. without any spectral analysis.

Despite the fact that the eigenfunctions of the Dirac-Coulomb operator
for $\nu=1$ can be explicitely constructed, nobody, to our knowledge,
has ever succeeded in deducing selfadjointness of this operator from that fact.

The main result of our paper states that for any potential $V$ for which an inequality like \eqref{R5} (with $c(V) \in (-1,1)$) holds true, we can define a distinguished self-adjoint extension of the operator $H_0+V$ defined on $C^\infty_c (\R^3, \C^4)$. In particular our analysis can be applied to the critical case $\nu=1$ of the Dirac-Coulomb Hamiltonian.
This extension will not necessarily satisfy condition \eqref{domain}, that is, the domain of $H_0+V$ will not necessarily be a space contained in $H^{1/2}(\R^3)$. For instance, this will be the case for the potential $V=-\frac1{|x|}$. The total energy, however, will still be finite for 
all functions in our domain.

Our method is not limited to Dirac operators with scalar potentials. For instance a Hardy-Dirac inequality for an hydrogenic atom in a constant magnetic field has been proved recently in \cite{Dolbeault-Esteban-Loss-06} provided the magnetic field is not too intense. As a consequence,  
a distinguished selfadjoint extension can be obtained also in this case.

\section{Main result and proof.}

In what follows we use the following assumption on the potential $V$.

\medskip
\noindent{\bf Assumption (A) :} {\sl 
$V:\R^3\to \R$ is a function such that for some constant $ c(V)\in (-1,1)$,  $\Gamma := \sup_{\R^3}V<1+c(V)$ and for every $\varphi \in C^\infty_c (\R^3, \C^2)$ ,} 
\be\label{R7} 
\int_{\R^3}\left(\frac{|{\boldsymbol{\sigma}}\cdot{\boldsymbol{\nabla}}\varphi|^2}{1+c(V)-V}\ +\ \big(1-c(V)+V\big)\,|\varphi|^2\right)\,dx\ \geq\ 0\; .
\ee

\medskip
In what follows $\gamma$ is any number in $(\Gamma, 1+c(V))$.
Consider the quadratic form
\be\label{beegamma}
b_\gamma(\varphi,\varphi) := \int_{\R^3}\left(\frac{|{\boldsymbol{\sigma}}\cdot{\boldsymbol{\nabla}}\varphi|^2}{\gamma - V}\ +\ \big(2-\gamma +V\big)\,|\varphi|^2\right)\,dx
\ee
defined on $C^\infty_c (\R^3, \C^2)$. Note that by assumption
\eqref{R7} this quadratic form is nonnegative and symmetric
on $C^\infty_c (\R^3, \C^2)$. Therefore it is closable and we denote its closure by $\widehat{b}_\gamma$ and its form domain by
$\mathcal{H}^\gamma_{+1}$. We denote by $S_\gamma$ the unique selfadjoint operator associated with $\widehat{b}_\gamma$: for
all $\varphi \in D(S_\gamma) \subset \mathcal{H}^\gamma_{+1}$,
\be
\widehat{b}_\gamma(\varphi, \varphi) = (\varphi, S_\gamma \varphi) \ .
\ee
$S_\gamma$ is an isometric isomorphism
from ${\mathcal H}^\gamma_{+1}$ to its dual ${\mathcal H}^\gamma_{-1}$. Using the second representation theorem in \cite{Kato} Theorem
2.23 we have that
${\mathcal H}^\gamma_{+1}$ is the operator domain of $S_\gamma^{1/2}$,
and
\be
\widehat{b}_\gamma(\varphi, \varphi) = (S_\gamma^{1/2}\varphi, S_\gamma^{1/2} \varphi) \ ,
\ee
for all $\varphi \in  {\mathcal H}^\gamma_{+1}$.
We now show that the above construction of ${\mathcal H}^\gamma_{+1}$
does not depend on $\gamma$. This follows from

\begin{prop} \label{prop-equivalence}
Under assumption {\bf (A)},
\be
\widehat{b}_\gamma(\varphi, \varphi) \le \widehat{b}_{\gamma'}(\varphi, \varphi)
+ [\gamma -\gamma']_+ \left[\frac{1} {(\gamma' -\Gamma)(\gamma -\Gamma)}+1 \right]\Vert \varphi \, \Vert^2_2
\ee
for all
$\gamma  , \ \gamma'$ in
$(\Gamma, 1+c(V))$ where $[\gamma -\gamma']_+ = \max \{\gamma -\gamma',0\} $. As a consequence, the spaces ${\mathcal H}^\gamma_{\pm 1}$ are independent of $\gamma$ and we denote them 
by ${\mathcal H}_{\pm 1}$.
\end{prop}

It suffices to prove the inequality for spinors $\varphi \in C^\infty_c (\R^3, \C^2)$. The proposition follows immediately from the elementary, pointwise inequality
\be
\frac{1}{\gamma -V(x)} - \frac{1}{\gamma' - V(x)}
\le \frac{[\gamma -\gamma']_+}{(\gamma - \Gamma)(\gamma' - \Gamma)} \ .
\ee
After this preparatory result we are ready to define an extension
of the Dirac operator originally defined on $C^\infty_c (\R^3, \C^4)$.

\begin{defi} \label{defi-domain}
The domain $\mathcal{D}$ of the Dirac operator  is the collection of all pairs
$\varphi \in {\mathcal H}_{+1}$, $\chi \in L^2(\R^3, \C^2)$ such that
\be \label{distributional}
(2-\gamma +V) \varphi -i\sigma\cdot\nabla\chi \ , \quad   -i\sigma\cdot\nabla \varphi + (V-\gamma)\chi\  \in \ L^2(\R^3, \C^2) \ . 
\ee
 
\end{defi}
The meaning of these two expressions is in the weak sense, i.e., the linear functional $\,(\eta, (V-\gamma)\chi) +(-i\sigma\cdot\nabla\eta, \varphi)\,$, which is defined for all test functions, extends uniquely to a bounded linear on $L^2(\R^3, \C^2)$. Likewise the same for 
$\,(\eta, (2-\gamma +V) \varphi) +(-i\sigma\cdot\nabla\eta, \chi)\,$.
From this definition it is clear that the domain does not depend 
on $\gamma$. Thus on this domain $\mathcal{D}$, we define the
Dirac operator as
\be
(H_0+V)\left(\begin{array}{c} \varphi \\ \chi \end{array} \right)= \left(\begin{array}{c}(V+1)\varphi -i\sigma\cdot\nabla\chi \\ -i\sigma\cdot\nabla\varphi+(V-1)\chi \end{array} \right) 
\ee
Note that for all vectors $(\varphi, \chi) \in \mathcal{D}$ the
expected total energy is finite. 

Our main result is the following theorem.
\begin{thm}\label{main}
Under the assumption {\bf (A)} on the potential $V$, the operator $H$, 
that is $H_0+V$ defined on $\mathcal{D}$ is selfadjoint. It is the unique selfadjoint extension of $H_0+V$ on $C^\infty_c (\R^3, \C^4)$  such that the domain is contained in $\mathcal{H}_{+1} \times L^2(\R^3, \C^2)$.  
\end{thm}

If in addition to {\bf (A)} we assume that \eqref{sotto1}
holds, it is natural to ask how our results relate to the ones
in \cite{Wust-73, Wust-75, Wust-77} and \cite{Nenciu-76}.
For simplicity we shall assume that $V(x) \le 0$ and that
\be \label{klauswuest}
\sup_{\R^3 \setminus \{0\}}|xV(x)| < 1 \ .
\ee
It was shown in \cite{Klaus-Wust-78} that for such potentials
$V(x)$ 
all these various selfadjoint extensions, which we call
$H_W, H_N$ coincide, i.e, $H_W=H_N=:H_{W,N}$ and the domain
is given by
\be
\mathcal{D}_{W}:=D(T^*) \cap D(r^{-1/2}) \ .
\ee
Here $T^*$ is the adjoint of $H_0+V$ restricted to $C^\infty_c(\R^3 \setminus \{0\}, \C^4)$. If we denote by $S^*$
the adjoint of $H_0+V$ restricted to $C^\infty_c(\R^3, \C^4)$,
then $S^* \subset T^*$. Note that $S^*$ is actually the operator $H_0+V$ defined on the 
set of spinors $(\varphi, \chi) \in L^2(\R^3, \C^2)$ for which \eqref{distributional} holds. 
Because of Theorem \ref{TR5} we find that ${\bf (A)}$ is satisfied
for such potentials. Hence we have a selfadjoint operator $H$ with
domain $\mathcal{D}$ which satisfies
\be
T \subset S \subset H=H^* \subset S^* \subset T^* \ ,
\ee 
moreover we have
\begin{thm}\label{notmain}
Assume that $V(x) \le 0$ and satisfies \eqref{klauswuest}.
Then $H=H_{W,N}$.
\end{thm}

The following two results are important in the proof of Theorem \ref{main}.

\begin{prop}\label{prop-scale}
Under assumption {\bf (A)} on the potential $V$, 
 \be\label{domain-estimate}
{\mathcal H}_{+1}\subset\Big\{\varphi\in L^2(\R^3, \C^2)\,:\; \frac{-i\sigma\cdot\nabla\varphi}{\gamma-V}\in L^2(\R^3,  \C^2)\Big\}\,,
\ee
where $ \nabla \varphi$ denotes the distributional gradient of $\varphi$.
Therefore, we have the `scale of spaces' \; ${\mathcal H}_{+1} \subset
L^2(\R^3, \C^2) \subset {\mathcal H}_{-1}$.
 \end{prop}
 
 \proof
 It suffices to prove the proposition for $\varphi \in C^\infty_c (\R^3, \C^2)$.
First we note that for all $x \in \R^3$
\be \label{in1}
\frac{1}{\gamma -V(x)} -\frac{1}{1+c(V) -V(x)} \ge   \frac{\delta}{(\gamma -V(x))^2}
\ee
where
\be
\delta := \frac{(\gamma -\Gamma)(1+c(V)-\gamma)}{1+c(V)-\Gamma} \ .
\ee
This follows from the assumption that $V(x) \le \Gamma$ for all $x$.
Noting that $1+c(V)-\gamma > \delta$, we get 
\be 
b_\gamma (\varphi, \varphi)\geq \delta\int_{\R^3}|\varphi|^2\,dx + \delta\int_{\R^3}\frac{|\sigma\cdot\nabla\varphi|^2}{(\gamma-V)^2}
\, dx\ .
\ee 
\finprf

\begin{lem} \label{lem-noidea}
For any $F$ in $L^2(\R^3, \C^2)$ and for any $\gamma \in (\Gamma, 1+c(V))$, 
\be
-i\sigma\cdot\nabla\left(\frac{ F}{\gamma-V}\right) \in \mathcal{H}_{-1} \ ,
\ee
where, once again, the gradient is to be interpreted in the distributional sense.
\end{lem}

\proof
By the definition of the distributional derivative, for every
$\eta \in C^\infty_c (\R^3, \C^2)$,
\be
\left|\left(-i\sigma\cdot\nabla \eta, \frac{ F}{\gamma-V} \right)\right|
= \left|\left(\frac{-i\sigma\cdot\nabla \eta}{\gamma-V}, F \right)\right| \le \Vert \eta \Vert_{\mathcal{H}_{+1}} \, \Vert F\Vert_2 \ .
\ee
Hence, the linear functional
\be
\eta \to \left(-i\sigma\cdot\nabla \eta, \frac{ F}{\gamma-V} \right)
\ee
extends uniquely to a bounded linear functional on $\mathcal{H}_{+1}$.
\finprf

\noindent{\sl Proof of Theorem \ref{main}.}
We shall prove Theorem \ref{main} by showing that for $\gamma \in (\Gamma, 1+c(V))$, $H+1-\gamma$ is symmetric and  a bijection
from its domain $\mathcal{D}$ onto $L^2(\R^3, \C^4)$. 
To prove the symmetry
we have to show that for both pairs $(\varphi, \chi)$, $(\tilde{\varphi}, \tilde{\chi})$ in the domain,
\be
\left((H+1-\gamma)\left(\begin{array}{c} \varphi \\ \chi \end{array} \right)\ ,\ 
\left(\begin{array}{c}\tilde{ \varphi} \\ \tilde{\chi} \end{array} \right)\right) =
((V-\gamma)\chi -i\sigma\cdot\nabla \varphi, \tilde{\chi}) + ((2-\gamma +V) \varphi -i\sigma\cdot\nabla \chi,
\tilde{\varphi})
\ee
equals
\be
(\chi, (V-\gamma)\tilde{\chi} -i\sigma\cdot\nabla \tilde{\varphi}) + (\varphi, (2-\gamma +V) \tilde{\varphi} -i\sigma\cdot\nabla \tilde{\chi}) \ =
\left(\left(\begin{array}{c} \varphi \\ \chi \end{array} \right)\ ,\ 
(H+1-\gamma)\left(\begin{array}{c}\tilde{ \varphi} \\ \tilde{\chi} \end{array} \right)\right)
.
\ee
First, note that since $(\varphi,\chi)$ is in the domain,
\be \label{elltwo}
(V-\gamma) \left[\chi +\frac{-i\sigma\cdot\nabla \varphi}{V-\gamma}\right] \in  L^2(\R^3, \C^2) \ .
\ee 
We now claim that
\be
((2-\gamma +V) \varphi -i\sigma\cdot\nabla \chi, \tilde{\varphi})
=(S_\gamma \varphi, \tilde{\varphi}) + \left((V-\gamma) \left[\chi +\frac{-i\sigma\cdot\nabla \varphi}{ V-\gamma} \right] ,\frac{-i\sigma\cdot\nabla \tilde{\varphi}}{V-\gamma}\right)
\ee
Note that each term makes sense. The one on the left, by definition
of the domain and the first on the right, because both $\varphi, \tilde{\varphi}$ are in ${\mathcal H}_{+1}$. The second term on the right
side makes sense because of \eqref{elltwo} above and Proposition \ref{prop-scale}.
Moreover both sides coincide for $\tilde{\varphi}$ chosen to be a test function and both are continuous in $\tilde{\varphi}$ with 
respect to the ${\mathcal H}_{+1}$ -norm. Hence the two expressions coincide on the domain. Thus we get
that
\be
\left((H+1-\gamma)\left(\begin{array}{c} \varphi \\ \chi \end{array} \right)\ ,\ 
\left(\begin{array}{c}\tilde{ \varphi} \\ \tilde{\chi} \end{array} \right)\right)
\ee
equals
\be
(S_\gamma \varphi, \tilde{\varphi}) + \left((V-\gamma) \left[\chi +\frac{-i\sigma\cdot\nabla \varphi}{V-\gamma} \right] ,\left[\tilde{\chi}+
\frac{-i\sigma\cdot\nabla \tilde{\varphi}}{ V-\gamma}\right]\right)
\ee
which is symmetric in $(\varphi,\chi)$ and $(\tilde{\varphi}, \tilde{\chi})$.
To show that the operator is onto,  pick any $F_1, F_2$ in $L^2(\R^3, \C^2)$. Since $S_\gamma$ is an isomorphism, there exists
a unique
$\varphi$ in ${\mathcal H}_{+1}$ such that
\be \label{phiequation}
S_\gamma \varphi = F_1 -i\sigma\cdot\nabla \left(\frac{F_2}{ \gamma -V}\right) \ .
\ee
Indeed, $F_1$ is in $L^2(\R^3, \C^2)$ and therefore in ${\mathcal H}_{-1}$. Moreover the second term is also in ${\mathcal H}_{-1}$ by Lemma \ref{lem-noidea}. 

Now define $\chi$ by
\be \label{chiequation}
\chi = F_2+ \frac{-i\sigma\cdot\nabla \varphi}{ \gamma - V}
\ee
which by Proposition \ref{prop-scale} is in $L^2(\R^3, \C^2)$. 

Now for any test function $\eta$ we have that
\be
(\eta, (V+2-\gamma)\varphi) + (-i\sigma\cdot\nabla\eta, \chi) =
(-i\sigma\cdot\nabla\eta, \frac{-i\sigma\cdot\nabla \varphi}{\gamma -V})+(\eta, (V+2-\gamma)\varphi) + (-i\sigma\cdot\nabla\eta, (\chi + \frac{-i\sigma\cdot\nabla\varphi}{V-\gamma}))
\ee
which equals
\be
(\eta,S_\gamma\varphi) + (-i\sigma\cdot\nabla \eta, (\chi + \frac{-i\sigma\cdot\nabla \varphi}{ V-\gamma}))
= (\eta, F_1)
\ee
This holds for all test functions $\eta$, but since $F_1$ is in $L^2(\R^3, \C^2)$,
the functional
$\eta \to (\eta, (V+2-\gamma)\varphi) + (-i\sigma\cdot\nabla \eta, \chi)$
extends uniquely to a linear continuous functional on $L^2(\R^3, \C^2)$ which
implies that
\be
(V+2-\gamma)\varphi -i\sigma\cdot\nabla \chi = F_1 \ .
\ee
Hence $(\varphi, \chi)$ is in the domain $\mathcal{D}$ and the operator $(H+1-\gamma)$ applied to $(\varphi, \chi)$ yields $(F_1, F_2)$.

Let us now prove the injectivity of $H+1-\gamma$. Assuming that
\be\label{esc3}
 \big(H+1-\gamma\big)\left(\begin{array}{cc}\varphi\\ \chi\end{array}\right)=\left(\begin{array}{cc}0\\0\end{array}\right)\,,
\ee
we find by \eqref {phiequation} and \eqref {chiequation}, 
$$\chi=\frac{-i\sigma\cdot\nabla\varphi}{\gamma-V}\;,\qquad S_\gamma \varphi=0\,.
$$
Since $S_\gamma$ is an isomorphism, this implies that $\varphi=\chi=0$.

It remains to show the uniqueness part in our theorem. Assume that
$C^\infty_c(\R^3,\C^4) \subset \mathcal{D}'$ is another selfadjoint extension such that whenever $(\varphi, \chi) \in \mathcal{D}'$,
the spinor $\varphi \in \mathcal{H}_+$ and, of course, $\chi \in L^2(\R^3,\C^2)$. Since $H$ is selfadjoint on this domain,
\be
\left(\left(\begin{array}{cc}\tilde{\varphi}\\ \tilde{\chi}\end{array}\right), H \left(\begin{array}{cc}\varphi\\ \chi\end{array}\right)\right)
= \left(H\left(\begin{array}{cc}\tilde{\varphi}\\ \tilde{\chi}\end{array}\right), \left(\begin{array}{cc}{\varphi}\\ {\chi}\end{array}\right)\right)
\ee 
for all $(\tilde{\varphi}, \tilde{\chi}) \in C^\infty_c(\R^3,\C^4)$.
This, however, means that the expressions \eqref{distributional}, defined in the distributional sense, belong to $L^2(\R^3,\C^2)$. 
Thus, $\mathcal{D}' \subset \mathcal{D}$ and hence
$\mathcal{D}' = \mathcal{D}$.
\finprf

\noindent
{\it Proof of Theorem \ref{notmain}}

If we can show that $\mathcal{D} \subset \mathcal{D}_{W}$
it follows that $ H \subset H_{W,N}$ and hence $ H = H_{W,N}$.
This follows from the next proposition in which we show that $\mathcal{D} \subset  D(r^{-1/2})$.

\begin{prop} \label{prop-equivalent}
Assume that $V(x) \le 0$ and satisfies \eqref{klauswuest}.
Then in the above notation $\mathcal{D} \subset  D(r^{-1/2})$.
\end{prop}
\noindent
{\it Proof}

We have by the assumption on $V(x)$, for any $\varphi \in C^\infty_c(\R^3,\C^2)$, 
\be
-\int V |\varphi|^2 dx \le \nu \int \frac{|\varphi|^2}{|x|} dx
\le \nu \int_{\R^3}\left(\frac{|{\boldsymbol{\sigma}}\cdot{\boldsymbol{\nabla}}\varphi|^2}{1+\frac{1}{|x|}}\ +\ \,|\varphi|^2\right)\,dx
\ee
where we have used \eqref{HD}. For $\gamma$ small enough we have
the elementary inequality
\be
\frac{\nu}{1+\frac{1}{|x|}} \le \frac{\nu^2}{\gamma -V} \ .
\ee
Hence
\begin{eqnarray}
&-&\nu \int \frac{|\varphi|^2}{|x|} dx  \le \nu ^2 \int_{\R^3}\frac{|{\boldsymbol{\sigma}}\cdot{\boldsymbol{\nabla}}\varphi|^2}{\gamma -V}\ +\ \,\nu \int_{\R^3} |\varphi|^2 \,dx \\
&=& \nu^2 b_\gamma(\varphi, \varphi) + (\nu+\nu^2(\gamma -2)) \int|\varphi|^2 dx - \nu^2 \int V|\varphi|^2 dx
\end{eqnarray}
Hence
\be
\nu \int \frac{|\varphi|^2}{|x|} dx \le \nu^2 b_\gamma(\varphi, \varphi) + (\nu+\nu^2(\gamma -2)) \int|\varphi|^2 dx + \nu^3 \int \frac{|\varphi|^2}{|x|} dx
\ee
which yields
\be
(1-\nu^2) \int \frac{|\varphi|^2}{|x|} dx \le \nu b_\gamma(\varphi, \varphi) + (1+\nu(\gamma -2)) \int|\varphi|^2 dx
\ee
This inequality extends to all $\varphi \in \mathcal{H}_{+1}$.
It remains to show a similar bound for $\chi$.

In the proof of the symmetry of $H$ we derived the formula
\be
\left(H\left(\begin{array}{c} \varphi \\ \chi \end{array} \right)\ ,\ 
\left(\begin{array}{c}{ \varphi} \\ {\chi} \end{array} \right)\right)
=(S_\gamma \varphi, {\varphi}) + \left((V-\gamma) \left[\chi +\frac{-i\sigma\cdot\nabla \varphi}{V-\gamma} \right] ,\left[{\chi}+
\frac{-i\sigma\cdot\nabla {\varphi}}{ V-\gamma}\right]\right)
\ee
which for $(\varphi, \chi) \in \mathcal{D}$ implies that
\be
0 \le \left(\left[\chi -\frac{-i\sigma\cdot\nabla \varphi}{\gamma-V} \right] ,(\gamma -V)\left[{\chi}-
\frac{-i\sigma\cdot\nabla {\varphi}}{\gamma -V}\right]\right) < \infty\ .
\ee
This implies, via simple estimates, that
\be
\left|\int_{\R^3} (\gamma-V)|\chi|^2 dx - \int_{\R^3} \frac{|\sigma\cdot\nabla \varphi|^2}{\gamma - V} dx\right| < \infty \ .
\ee
Since $\varphi \in \mathcal{H}_{+1}$, $(\varphi, S_\gamma \varphi)$
is finite and since 
$\int \frac{|\varphi|^2}{|x|} dx <\infty$, we also have that
\be
\int_{\R^3} \frac{|\sigma\cdot\nabla \varphi|^2}{\gamma - V} dx <\infty \ .
\ee
This implies the result.
\finprf 
\bigskip\noindent{\bf Acknowledgments.} M.J.E. would like to thank M. Lewin,  E. S\'er\'e and J.-P. Solovej for various discussions on the self-adjointness of Dirac operators.

M.J.E. and M.L. wish to express their gratitute to  Georgia Tech and Ceremade for their hospitality. M.J.E. acknowledges support from ANR Accquarel project and European Program ``Analysis and Quantum'' HPRN-CT \# 2002-00277. M.L. is partially supported by U.S. National Science Foundation grant DMS DMS 06-00037.

\medskip\noindent\copyright\, 2007 by the authors. This paper may be reproduced, in its entirety, for non-commercial purposes.


\end{document}